\documentclass[a4paper]{amsart}

\usepackage{microtype}
\usepackage{xspace}

\usepackage{hyperref,url}
\usepackage{listings}

\usepackage{graphicx}
\usepackage{tikz}
\usetikzlibrary{shapes,calc,arrows,through,intersections}

\newcommand\RR{{\mathbb R}}
\newcommand\MatchTheNet{\texttt{MatchTheNet}\xspace}
\newcommand\polymake{\texttt{polymake}\xspace}
\lstdefinestyle{polymake}{
language={Perl},
xleftmargin=2em,
xrightmargin=2em,
basicstyle=\linespread{1.3}\ttfamily,
identifierstyle=,
keywordstyle=\bfseries,
commentstyle=\color{gray},
moredelim=[is][\color{gray}]{||}{||},
stringstyle=\ttfamily,
breaklines=true,
sensitive=true,
columns=fullflexible,
numbers=none, numberstyle=\tiny, frame=none, fontadjust=false, tabsize=4,
morekeywords={},
otherkeywords={polytope>,ideal>},
keywordstyle={\color{black!80}\bfseries},
string=[b]{"},
morestring=[b]{'},
showstringspaces=false,
stringstyle={\color{gray}},
classoffset=2,
morekeywords={archimedean_solid, fan, planar_net, threejs, min, polytope>},
keywordstyle={\color{black!70}\bfseries},
classoffset=3,
morekeywords={VISUAL, MAXIMAL_POLYTOPES},
keywordstyle={\color{black!50}\bfseries},
classoffset=4,
morekeywords={},
keywordstyle={\color{black!70}\bfseries\slshape},
classoffset=5,
morekeywords={VertexLabels, VertexColor,FacetTransparency,FacetColor},
keywordstyle={\color{black!90}\slshape}
}

\title{\MatchTheNet\ - an educational game on 3-dimensional polytopes}
\thanks{Research by M. Joswig is supported by Einstein Foundation Berlin and Deutsche Forschungsgemeinschaft (Priority Program 1489: ``Experimental methods in algebra, geometry, and number theory'',
  SFB/TRR 109: ``Discretization in Geometry and Dynamics'' and SFB/TRR 195: ``Symbolic Tools in Mathematics and their Application'')}

\author{Michael Joswig}
\author{Georg Loho}
\author{Benjamin Lorenz}
\author{Rico Raber}
\address{Institut f\"ur Mathematik, MA 6-2, Technische Universit\"at Berlin, \\
  Straße des 17.Juni 136, 10623 Berlin, Germany \\
  \texttt{\{joswig,loho,lorenz,raber\}@math.tu-berlin.de}
}


\begin{document}

\maketitle

\begin{abstract}
  We present an interactive game which challenges a single player to match 3-dimensional polytopes to their planar nets.
  It is open source, and it runs in standard web browsers.
\end{abstract}

\section{Introduction}
A \emph{polytope} is the convex hull of finitely many points in Euclidean space.
While their study goes back to antiquity, polytopes are still an active research topic; see, e.g., Ziegler~\cite{Ziegler95}.
The \emph{dimension} of a polytope is the dimension of its affine span.
The first non-trivial class of polytopes are the $3$-polytopes, i.e, those of dimension three.
This includes the Platonic and Archimedean solids as their most prominent examples.
The combinatorics of a $3$-polytope $P$ is determined by its (vertex--edge) graph $\Gamma$, which is planar.
Our game \MatchTheNet invites to play with these geometric objects.
It is based on the infrastructure of the software system \polymake \cite{DMV:polymake}.

\begin{figure}[htb]
  \centering
  \begin{minipage}[b]{0.3\textwidth}
    \includegraphics[width=\textwidth]{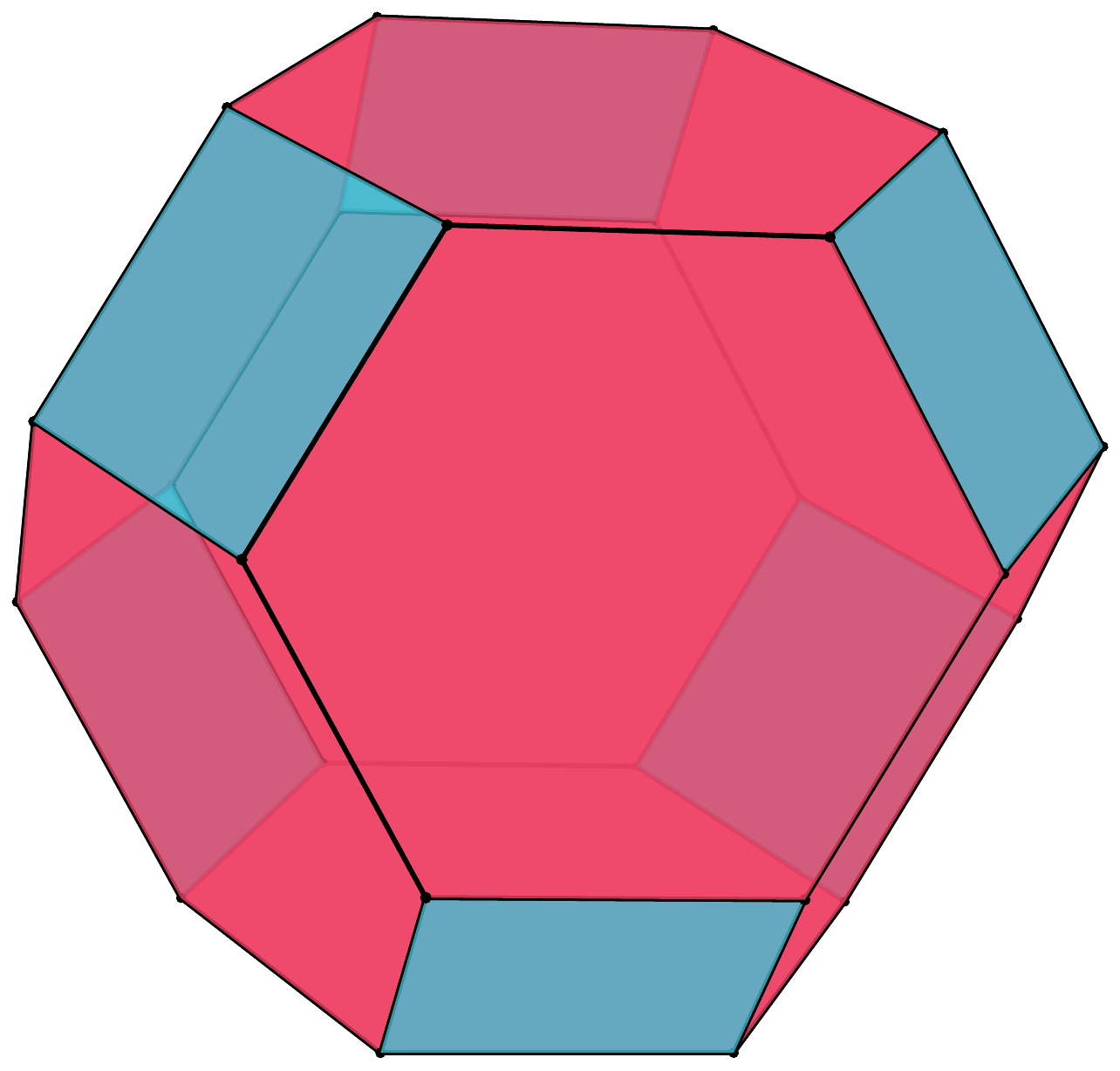}
  \end{minipage}
  \begin{minipage}[b]{0.3\textwidth}
    \includegraphics[width=\textwidth]{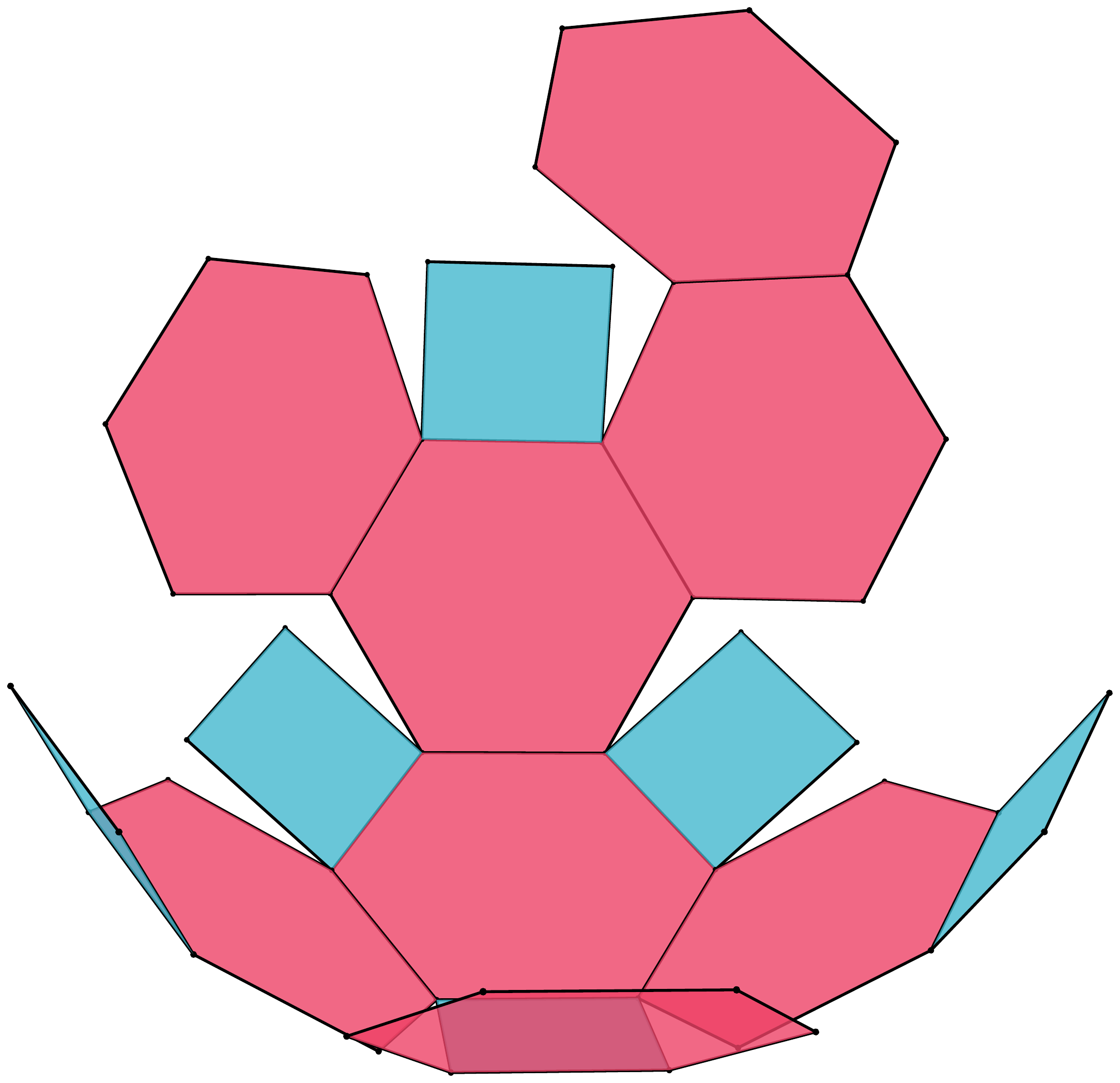}
  \end{minipage}
  \begin{minipage}[b]{0.3\textwidth}
    \includegraphics[width=\textwidth]{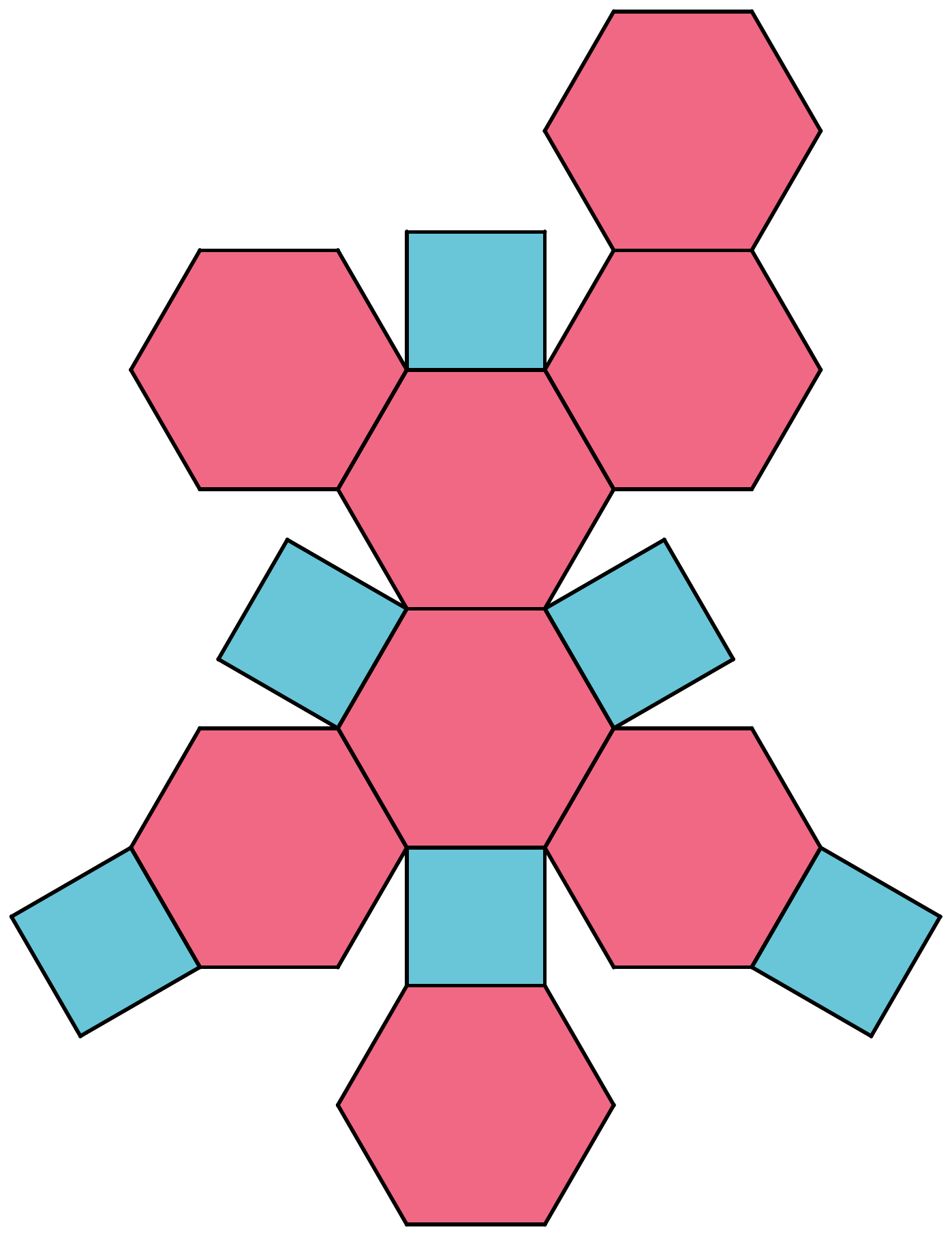}
  \end{minipage}
  \caption{Unfolding the truncated octahedron, which is an Archimedean solid.}
  \label{fig:unfolding}
\end{figure}

We obtain a \emph{planar net} of a $3$-polytope by cutting the boundary along several edges. The resulting shape unfolds to a flat and connected figure, similar to Fig.~\ref{fig:unfolding}. 
Given a planar net on a sheet of paper, one can cut out the shape, fold it along sketched edges and glue it along some boundary edges to regain the original polytope.
In the following, planar nets are described more formally.

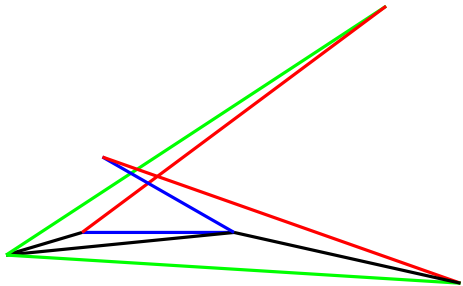
\begin{figure}[bht]\centering
\begin{tikzpicture}

\tikzset{edgesty/.style={black, very thick}}

\coordinate (P1) at (0,0);
\coordinate (P2) at (-2,0);
\coordinate (P3) at (-3,-0.3);
\coordinate (P4) at (2,3);
\coordinate (P5) at (150:2);

\path[name path = c1] let \p1 = ($(P3) - (P4)$), \n1 = {veclen(\x1,\y1)} in (P3) ++(0:\n1) arc (0:-10:\n1);

\path[name path = c2] let \p2 = ($(P2) - (P4)$), \n2 = {veclen(\x2,\y2)} in (P5) ++(-15:\n2) arc (-15:-25:\n2);

\path[name intersections={of=c1 and c2,by={P6}}];

\draw[edgesty,blue] (P1) -- (P2);
\draw[edgesty] (P2) -- (P3);
\draw[edgesty] (P3) -- (P1);
\draw[edgesty,green] (P3) -- (P4);
\draw[edgesty,red] (P2) -- (P4);
\draw[edgesty,blue] (P1) -- (P5);
\draw[edgesty,red] (P5) -- (P6);
\draw[edgesty,green] (P3) -- (P6);
\draw[edgesty] (P1) -- (P6);

\end{tikzpicture}
\caption{Tetrahedron with an attempt to unfold that fails.
  It arises from a spanning tree in the dual graph which is a path.
  Each spanning tree which has a node of degree three gives a proper planar net; this works for any tetrahedron.}
\label{fig:tetrahedron}
\end{figure}

The \emph{dual graph} $\Delta$ of a $3$-polytope $P$ is the abstract graph which has the facets of $P$ as nodes, while its edges are given by those pairs of facets which share a common edge.
If we pick a spanning tree $T$ of $\Delta$, then, as in von Staudt's proof of Euler's formula, the edges of $\Gamma$ which are not dual to any edge in $T$ form a spanning tree $T^*$ of $\Gamma$.
In topological terms, the complementary pair $(T,T^*)$ of spanning trees corresponds to the two critical points of an optimal Morse function of the $2$-sphere.
We may view the dual spanning tree $T^*$ as a subset of the boundary $\partial P$.
Then we obtain a map from $\pi:\partial P\setminus T^* \to \RR^2$ as follows.
We start out by picking a facet $R$ of $P$ and map it isometrically into the plane.
Then, for each facet $F$ adjacent to $R$ there is a unique way of extending this map such that it is an isometry if restricted to $F$.
Now $\pi$ is defined inductively by following the unique path from any facet to the root facet $R$ in the spanning tree $T$.
The dual tree $T^*$ is said to define an \emph{edge cutting}, and the map $\pi$ only depends on $T^*$, but not on the choice of the root facet $R$.
If $\pi$ is injective, the closure of the image $\pi(\partial P\setminus T^*)$ is called a \emph{planar net} (or an \emph{unfolding}) of $P$.
It is an interesting open problem, whether or not each $3$-polytope admits a planar net; see, e.g., \cite{Shephard:1975} and \cite{BDEKMS:2003} as well as the monograph \cite{DemaineOrourke:2007} for an overview of topics related.
Figure~\ref{fig:tetrahedron} shows that an attempt to unfold may fail.

\MatchTheNet is a game where a single player is asked to match a set of planar nets to a set of $3$-polytopes.
The difficulty ranges from easy (suitable for kids in elementary school) to hard (recreational puzzle for grown-up mathematicians).
The game mechanics is written in \texttt{JavaScript}, and it runs in any web browser, either locally or over the Internet.
It can be played online at \url{www.matchthenet.de}, downloaded at \url{https://github.com/polymake/matchthenet}, and it is part of the \texttt{Imaginary} project.

\section{Playing the Game}

The front page of \MatchTheNet explains the rules, and the player can choose the language, the level of difficulty and the number of polytopes per round.
That number, which we will refer to as $k$ here, ranges between two and five.
One game lasts for five rounds.
In each round the player sees $k$ polytopes in the top row of the screen and $k$ planar nets in the bottom row.
The player swaps the planar nets with the mouse until she is confident that each polytope sits right above its planar net.
Hitting the ``submit'' button reveals the score, which is the total number of correct matches.
Afterwards the player can either look at the solution or continue with the next round.
After the fifth round the final score is displayed and compared to the current high score.
During the game the polytopes can be rotated freely with mouse to look at them from all sides.

There are various ways to make the game easier or more difficult.
We offer seven levels.
In general, the more facets the polytope has the more difficult it is to recognize.
Further, it makes a difference if the coloring of the facets gives some guidance to the combinatorics.
For instance, on Level 5 there are polytopes which come from a random construction, but color helps to identify the number of vertices on each facet.
On Level 6 the polytopes are the same, but all facets are green.
The highest Level 7 has triplets of polytopes chosen by hand, which are very similar to one another. For this level, only $k=3$ is available.

\section{Our Collection of Polytopes}

The bulk of our pre-computed $3$-polytopes are regular polytopes and their generalizations.
A \emph{Platonic solid} (or \emph{regular} $3$-polytope) admits an automorphism group (of rigid motions) which acts transitively on the set of maximal flags, i.e., the triplets consisting of a vertex, an edge and a facet which are incident; there are five combinatorial types.
More generally, the \emph{Johnson solids} are the $3$-polytopes whose facets are regular polygons of various gonalities.
An \emph{Archimedean solid} (or \emph{semi-regular} $3$-polytope) is a Johnson solid which admits a vertex-transitive group; there are 13 combinatorial types in addition to the regular ones.
The \emph{Catalan solids} are the duals of the Archimedean solids.
There are 92 combinatorial types of \emph{proper Johnson solids}, i.e., those which are not Archimedean \cite{Johnson66}.
Taking also the duals of the proper Johnson solids into account we arrive at five classes of $3$-polytopes which are pairwise disjoint.
Their numbers add up to $5+13+13+92+92=215$.
All of them are contained in the data base of \MatchTheNet.
Figure~\ref{fig:unfolding} shows an unfolding of an Archimedean solid.

Additionally, we computed fifty random $3$-polytopes by the following two-step procedure.
In the first step we choose hyperplanes tangent to the unit sphere uniformly at random.
Almost surely the resulting polytope $Q$ is \emph{simple}, i.e., each vertex is contained in precisely three facets.
In the second step we pick a certain portion of the vertices of $Q$, again uniformly at random, take their convex hull, and this is our random polytope.
Usually, such a polytope is neither simple nor dual to simple, i.e., simplicial.

For each level there is a subset of the entire collection from which polytopes are chosen at random during the game.
The highest level is different in that certain triplets of Johnson polytopes are chosen by hand.

\section{Computations in \polymake}

\polymake is open source software for research in polyhedral geometry \cite{DMV:polymake}.
It deals with polytopes, polyhedra and fans as well as simplicial complexes, matroids, graphs, tropical hypersurfaces, and other objects.
For \MatchTheNet we use \polymake as an engine to pre-compute all $3$-polytopes and their planar nets used in our game.

To give an example we show how to produce the planar net of the truncated octahedron shown in Figure~\ref{fig:unfolding} to the right.
This code is valid for \polymake version 3.0 or higher.
First we construct the polytope and its planar net.
The latter employs a heuristic with backtracking.

\begin{lstlisting}[style=polymake,aboveskip=1em,belowskip=1em,linerange=1-2]
polytope> $polytope = archimedean_solid('truncated_octahedron');
polytope> $net = fan::planar_net($polytope);
\end{lstlisting}

For visualization \polymake offers several backends.
Here, as for \MatchTheNet, we use the library \texttt{three.js} for a direct rendering in a web browser.

\begin{lstlisting}[style=polymake,aboveskip=1em,belowskip=1em,linerange=3-9]
polytope> @colors = ('green','blue','purple','red','grey');
polytope> threejs( $net->VISUAL( VertexLabels => "hidden", 
                 VertexColor => "black", 
                 FacetTransparency => 0.8, 
                 FacetColor => sub { 
                    $colors[ min($net->MAXIMAL_POLYTOPES->[shift]->size-3, @colors-1) ] 
                 } ));
\end{lstlisting}

In this example the color of each facet is determined by its number of vertices.
So triangles become green, quadrangles blue, pentagons purple and hexagons red; all others will be shown in gray.

\subsection*{Acknowledgements}
We are indebted to the \texttt{Imaginary} team for a lot of inspiration and fruitful discussions during the design of \MatchTheNet.

\bibliography{main}

\end{document}